\documentclass[onecollarge]{STJour2}        		% onecolumn
\smartqed 
\usepackage[square,sort&compress,comma,numbers]{natbib}	% fuer BibTeX
\usepackage{graphicx}
\usepackage{amssymb}
\usepackage{amsmath}
                             
\usepackage{multicol}        
\usepackage[bottom]{footmisc}

\usepackage{amsmath}
\usepackage{amsfonts}
\usepackage{url}
\usepackage[labelfont=bf,labelsep=space]{caption}
\usepackage{subcaption}
\usepackage{microtype}
\usepackage{mathtools}
\usepackage{booktabs}

\begin{document}
\unitlength1cm

\title{A Particle-based Multiscale Solver for Compressible Liquid-Vapor Flow}

\author{Jim Magiera \and Christian Rohde}
\institute{Jim Magiera, 
{jim.magiera@mathematik.uni-stuttgart.de}
\\
Christian Rohde, {crohde@mathematik.uni-stuttgart.de}
\\ University of Stuttgart, Institute for Applied Analysis and Numerical Simulation}

\date{April 2018}

\maketitle

\abstract{
To describe complex flow systems accurately, it is in many cases important to account for the properties of fluid flows on a microscopic scale.  
In this work, we focus on the description of liquid-vapor flow with a sharp interface between the phases.
The local phase dynamics at the interface can be interpreted as a Riemann problem for which we develop a multiscale solver in the spirit of the heterogeneous multiscale method (HMM) \cite{e.engquist.ea:heterogeneous:2007}, using a particle-based microscale model to augment the macroscopic two-phase flow system.
The application of a microscale model makes it possible to use the intrinsic properties of the fluid at the microscale, instead of formulating (ad-hoc) constitutive relations.
\keywords{Multiscale modeling, heterogeneous multiscale method, conservation laws, compressible two-phase flow, liquid-vapor flow, sharp interface resolution, Riemann problem, particle chain model, model reduction, machine learning.}
}

\section{Introduction}

For many problems in science and engineering microscopic properties can heavily influence the macroscopic behavior. 
Therefore it is important to consider microscopic effects in the mathematical model development. 
The obvious possibility to account for such small-scale effects is to solve the microscopic model everywhere. 
However, despite advances in computing power over the last decades, it is usually still not feasible. 
This scenario applies to the case of compressible fluid flows with liquid-vapor phase transition. 
Most applications require a computational domain on a laboratory scale, which is many orders of magnitude apart from a truly microscopic model that considers effects on the molecular level.
\newline
One approach to this problem is to perform multiscale domain-decomposition of micro- and macroscale models, where  in a part of the domain a microscale particle model is solved instead of the macroscale model, and both models are coupled via suitable  boundary conditions. 
This coupling approach has been investigated for example in \cite{ren:analytical:2007} for the incompressible Navier--Stokes equations on the macroscale and a Lennard--Jones particle model as the microscale model.
In \cite{li.yang.ea:multiscale:2010} multiscale domain-decomposition is applied for crack propagation in brittle materials, where a set of conservation laws is used in the continuum domain and near the crack a microscale particle model is applied.
Furthermore, the phase change of a liquid on a hot plate has been examined in \cite{cosden:hybrid:2013}. 
\newline
In this work however, we propose a multiscale model for the description of single liquid droplets, based on the heterogeneous multiscale method (HMM) 
\cite{e.engquist.ea:heterogeneous:2007,e:principles:2011}, which is a general framework for developing multiscale models. 
The main idea behind it is to compute solutions of a microscopic model for some given macroscopic constraints, and propagate hereby obtained parameters to the macroscopic model. 
Consequently, instead of performing multiscale domain-decomposition coupling of the scales, a data-based approach is promoted.

\section{The Macroscale Model: Compressible, Isothermal Euler Equations}

\begin{figure}
% \sidecaption
\centering
\includegraphics[width=0.2\columnwidth]{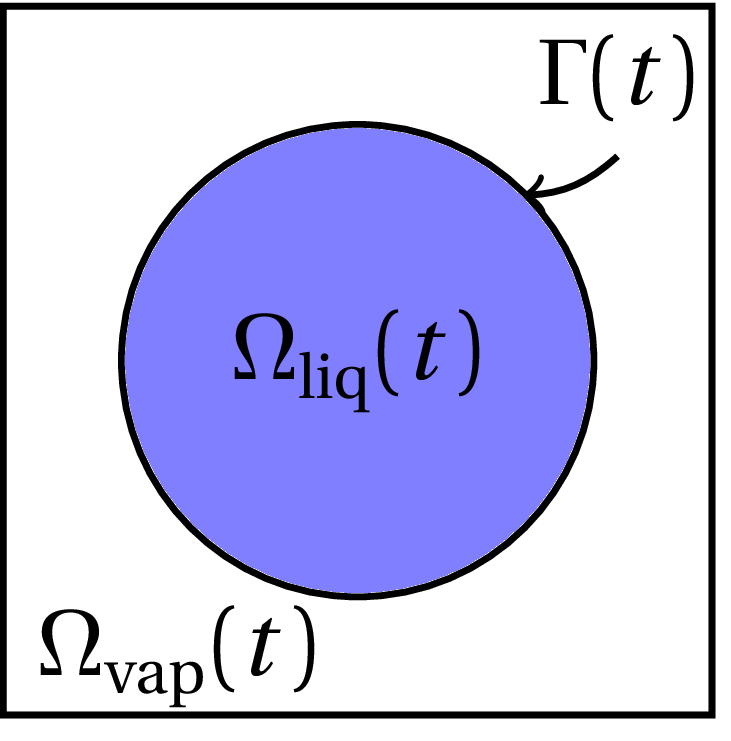}
\caption{Sketch of the two-phase flow domains.}
\label{fig:droplet_domains}
\end{figure}

On the macroscopic scale we consider the behavior of a single liquid droplet in a vapor atmosphere.
For such two-phase flows it is possible to consider either a diffuse interface approach \cite{anderson.mcfadden.ea:diffuse:1998}, where the phase-boundary has a finite thickness, or a sharp interface approach, as in \cite{rohde.zeiler:relaxation:2015,zeiler:liquid:2015}, where a discontinuous transition between the phases is present. 
In this work, we follow the second approach and assume that the interface between the phases is represented as a discontinuous shock wave.
\newline
Furthermore, we assume that the fluid flow is compressible, inviscid and isothermal at reference temperature $T_{\mathrm{ref}}$, such that the dynamics are described by the isothermal Euler equations
\begin{align} \label{eq:euler-2phase} 
\begin{aligned}
  \partial_t \rho + \nabla \cdot (\rho v) &= 0, \\
  \partial_t (\rho v) + \nabla \cdot (\rho v \otimes v) + \nabla p(\rho) &= 0,
  \end{aligned}
\end{align}
for the density $\rho$ and velocity $v$ in the space-time domain $\Omega \times (0,T)$, with $T > 0$ and $\Omega \subset \mathbb{R}^d$ an open set.
\newline
To describe the two separate phases, we distinguish at each point of time $t \in [0,T]$, between the two distinct bulk phases $\Omega_{\mathrm{vap}}(t)$ and $\Omega_{\mathrm{liq}}(t)$ with common boundary/interface $\Gamma(t)$, such that $\Omega_{\mathrm{vap}}(t) \cup \Omega_{\mathrm{liq}}(t) \cup \Gamma(t) = \Omega$. 
Figure \ref{fig:droplet_domains} shows a sketch of this setting. 
%\newline
To close the system \eqref{eq:euler-2phase}, the pressure $p$ has to be specified. 
For describing a generic two-phase system we consider the van der Waals pressure function, in terms of the specific volume $\tau = \tfrac{1}{\rho}$, as in \cite{rohde.zeiler:relaxation:2015},
\begin{align} \label{eq:van-der-waals-pressure}
p(\tau) =  \frac{R T_{\mathrm{ref}} }{\tau - b} - \frac{1}{\tau^2},
\end{align}
with some constants $R,b,a >0$. 
If the temperature $T_{\mathrm{ref}}$ is greater than the critical temperature 
$T_{\mathrm{c}} = \tfrac{8 a}{27 R b}$
the van der Waals pressure function is monotone and the system \eqref{eq:euler-2phase} is hyperbolic. 
However, if $T_{\mathrm{ref}} < T_{\mathrm{c}}$, the pressure is non-monotone, and the system becomes elliptic for $\tau \in (\tau^{\mathrm{max}}_{\mathrm{liq}},\tau^{\mathrm{min}}_{\mathrm{vap}})$ which is called spinodal region. 
Thus, we define the admissible set of densities as $\mathcal{A}_{\mathrm{vdw}} := (b,\infty) \setminus  (\tau^{\mathrm{max}}_{\mathrm{liq}},\tau^{\mathrm{min}}_{\mathrm{vap}})$ and distinguish between the liquid phase for $\tau \in (b, \tau^{\mathrm{max}}_{\mathrm{liq}})$ and the vapor phase for $\tau \in (\tau^{\mathrm{min}}_{\mathrm{vap}}, \infty)$.

\begin{figure}
% \sidecaption
\centering
\includegraphics[width=0.45\columnwidth]{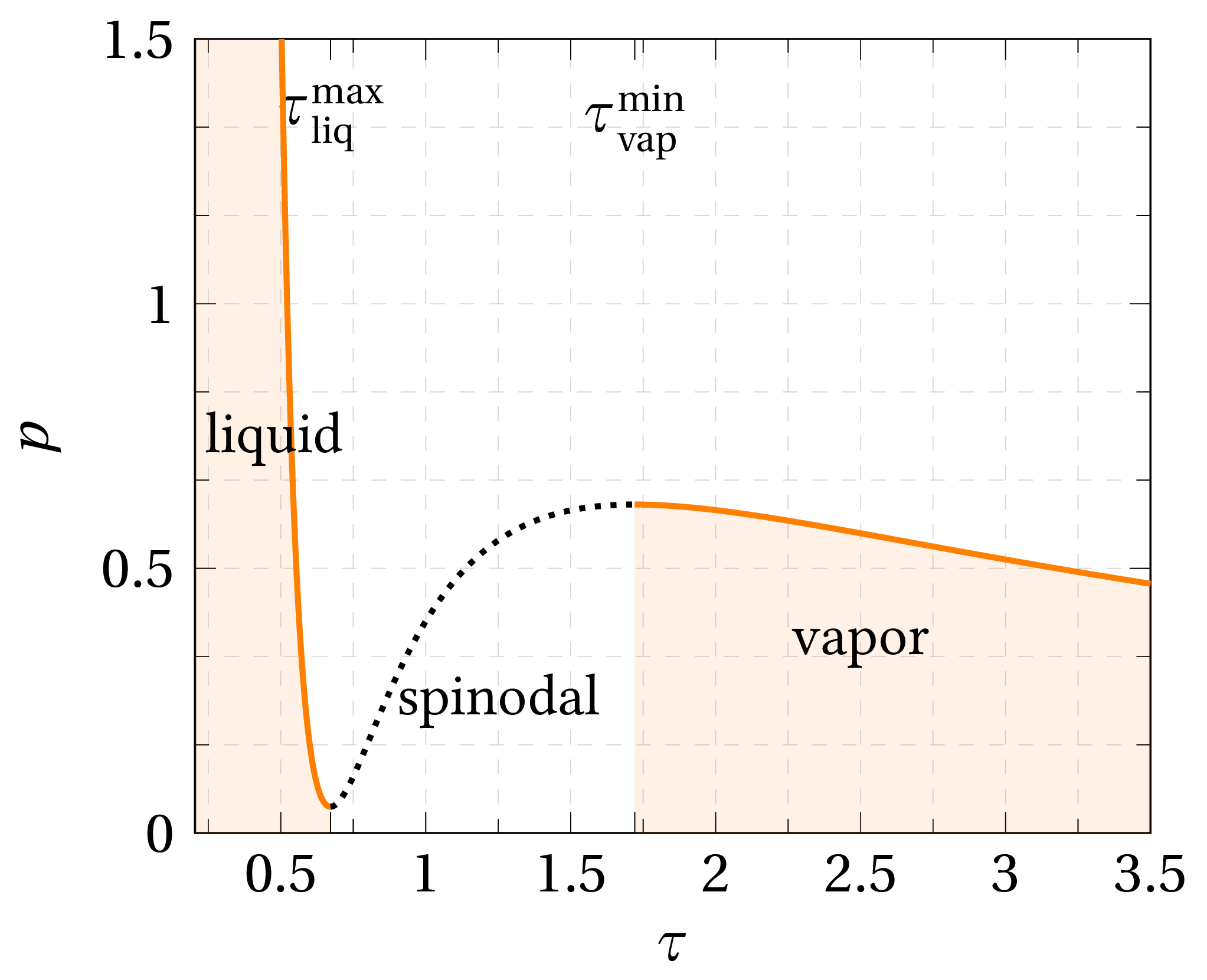}
\caption{The van der Waals pressure function for $T_{\mathrm{ref}} < T_{\mathrm{c}}$.}
\label{fig:vdw-pressure}
\end{figure}
 
\begin{figure}[h]
\sidecaption
\centering
\mbox{}\hfill
\includegraphics[width=0.3\columnwidth]{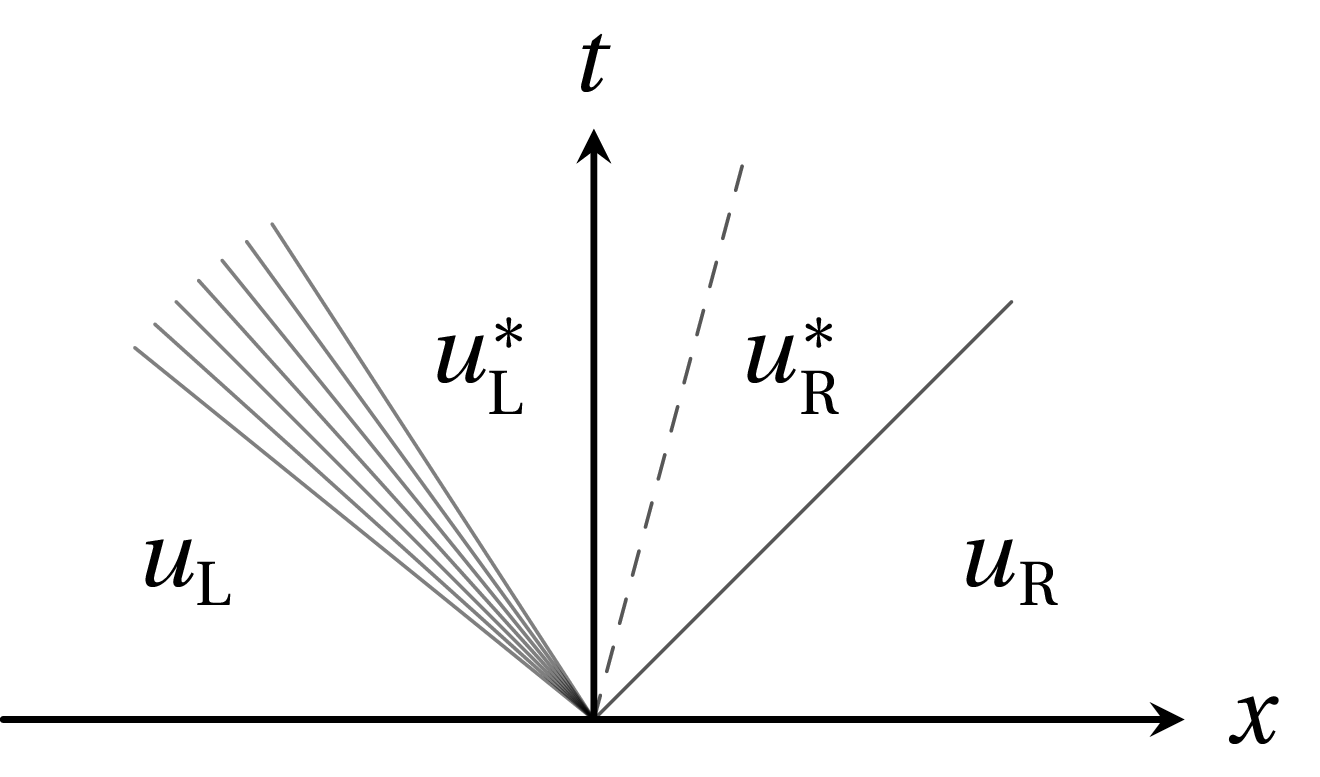}
\hfill\mbox{}
\caption{Sketch of a wave pattern for two-phase flow. The dashed line indicates the phase transition, which is sharp as an additional discontinuous wave.}
\label{fig:wavepattern}
\end{figure}

In order to complete the two-phase model, we have to formulate, besides initial and boundary conditions, some additional coupling conditions at the interface $\Gamma(t)$. 
Therefore, let $\xi \in \Gamma(t)$ and $t \in [0,T)$ be fixed. 
The speed of the interface $\Gamma(t)$ in normal direction $\vec{\nu}(\xi, t) \in \mathbb{S}^{d-1}$ (always pointing into the vapor phase) is denoted by $s(\xi,t) \in \mathbb{R}$. 
Then the mass and momentum balance at the interface, neglecting surface tension, take the following form
\begin{align}  \label{eq:interface_conditions}
\begin{aligned}
 \mbox{} [{\!}[ \rho ( v \cdot \vec{\nu} - s ) ]{\!}] &= 0, \\ 
 \mbox{} [{\!}[ \rho ( v \cdot \vec{\nu} - s ) v \cdot \vec{\nu} + p(\rho) ]{\!}] &= 0, \\ %(d-1) \kappa \gamma, \\ %
 \mbox{} [{\!}[ v \cdot \vec{t}]{\!}] &= 0, \quad \forall \,  \vec{t} \perp \vec{\nu},
 \end{aligned}
\end{align}
where $[{\!}[ \, \cdot \, ]{\!}]$ denotes the difference between liquid and vapor phase values.
The well-posedness of the free boundary value problem requires still another coupling condition. 
For the relevant subsonic case one assumes that this condition can be written down as an algebraic equation, called kinetic relation. 
It describes the entropy dissipation at the interface \cite{truskinovsky:kinks:1993}. 
\newline
For given initial Riemann data $u_{\mathrm{L}} = (\rho, \rho v)_{\mathrm{L}}$ for $x \leq 0$, and $u_{\mathrm{R}} = (\rho, \rho v)_{\mathrm{R}}$ for $x > 0$, the solution of the initial value problem \eqref{eq:euler-2phase} evolves (in contrast to the one-phase case) as a 3-wave pattern -- a sketch of such a wave pattern is depicted in Figure \ref{fig:wavepattern}.

Two-phase models with kinetic relations have been investigated in detail, see for example
\cite{abeyaratne.knowles:kinetic:1991,bedjaoui.chalons.ea:non:2005,merkle:sharp:2007}. 
\newline
However, it can be seen that for certain settings, the wrong choice of the kinetic relation can lead to a behavior of the model that is not observed by physical experiments, see e.g. \cite{zeiler:liquid:2015}. 
For that reason, we want to return to a more elementary notion of the physical properties and regard the flow at the interface on a molecular level. 
This has the advantage that no kinetic relation is needed. 
Furthermore, most physical parameters on the molecular level can be determined accurately by experiments. 
These advantages become even more apparent if one considers non-isothermal multiphase flow and mixtures, where the physically correct choice of the kinetic relation is usually not clear.

\section{The Microscale Model: Particle Chain Model}
\label{sec:particle_model}

\begin{figure}
\sidecaption
\centering
\mbox{}\hfill
\includegraphics[width=0.45\columnwidth]{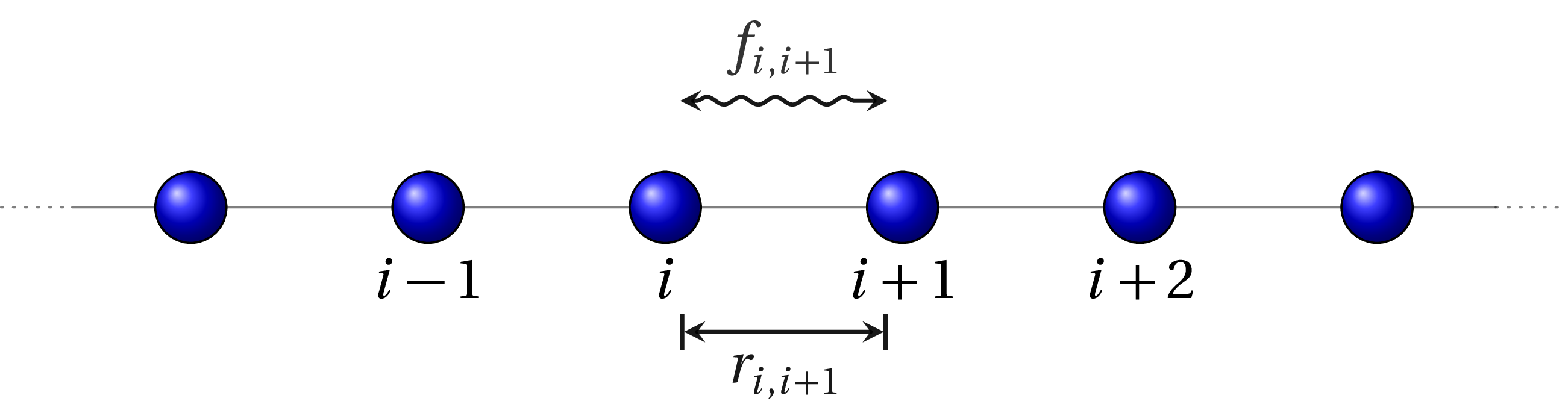}
\hfill\mbox{}
\caption{Sketch of the particle chain model.}
\label{fig:particlechain}
\end{figure}

For the description of the liquid-vapor interaction of droplets on a microscopic scale, we apply an atomistic one-dimensional particle chain model, which has been investigated for example in \cite{herrmann.rademacher:riemann:2010}.
More precisely that means that we consider a one-dimensional system of $N$ particles with position $x_i = x_i(t)$, velocity $v_i = v_i(t)$ and mass $m_i$, for $i=1,\ldots,N$. 
The distance between the $i$-th and $(i+1)$-th particle is given by $r_{i,i+1} = \lvert x_{i+1} - x_i \rvert$, see Figure \ref{fig:particlechain}.
The particles are assumed to interact only with direct neighbors via a potential $\phi \colon \mathbb{R}^+ \to \mathbb{R} : r \mapsto \phi(r)$, where $r$ denotes the distance between the particles.
The $i$-th particle is subject to the forces $f_{i-1,i}$, $f_{i,i+1}$ originating from the potentials of the neighboring particles, the resulting force $f_i$ is therefore given by
\begin{align*}
 f_i = f_{i-1,i} + f_{i,i+1} 
= \phi'(\lvert x_{i-1} - x_i \rvert) - \phi'(\lvert x_{i+1} - x_i \rvert ).
\end{align*}
Consequently, the acceleration $a_i = a_i(t)$ of the $i$-th particle is given by $a_i = f_i / m_i$. 
For the boundary conditions we assume that $f_0$ and $f_N$ are zero.
This gives us the following ordinary initial value problem for the particle motion
\begin{align} 
  \tfrac{\mathrm{d}^2}{\mathrm{d} t^2} x_i(t) &= \tfrac{1}{m_i} f_i(t), 
  &
  x_i(0) &= x^0_i, & v_i(0) &= v^0_i, 
\end{align}
with initial positions $x^0_i$ and velocities $v^0_i$ for $i = 1,\ldots,N$.

\subsection{Micro-/Macroscale Conversion: Irving--Kirkwood Formulas}

To design a multiscale scheme that accounts for microscopic properties it is essential to convert the key quantities from the macroscopic to the microscopic scale and vice versa. 
In case of a particle model this can be achieved via the  Irving--Kirkwood formulas \cite{irving.kirkwood:statistical:1950}. 
The microscopic instantaneous density $\rho(x,t)$ and momentum $(\rho v)(x,t)$ distributions are realized by 
\begin{align} \label{eq:irving-kirkwood:1}
\begin{aligned}
 \widehat{\rho}(x,t) & = \sum_{i=1}^{N} m_i \, \delta(x - x_i(t)), 
 &
 (\widehat{\rho} \widehat{v})(x,t) & = \sum_{i=1}^{N} m_i \, v_i(t) \, \delta(x - x_i(t)),     
 \end{aligned}
\end{align}
where $m_i$, $x_i$, $v_i$ are the mass, position and velocity of the $i$-th particle and $\delta$ denotes the Dirac distribution. 
The instantaneous pressure distribution $\widehat{p}(x,t)$ is given by 
\begin{align*}
 \widehat{p}(x,t) = \frac{1}{d} \biggl( \sum_{i=1}^{N} (m_i \overline{v}_i \cdot \overline{v}_i) \, \delta(x - x_i(t)) + \sum_{\substack{i = 1,\ldots,N \\  j < i}} (f_{ij} \cdot r_{ij}) \, \lambda_{ij}(x,t) \biggr),
\end{align*}
Here, $\overline{v}_i$ denotes the relative velocity with respect to a local mean value, $r_{ij}(t) := (x_i(t) - x_j(t))$,  and $\lambda_{ij}(x,t)$ is defined as
\begin{align*}
 \lambda_{ij}(x,t) := \int^{1}_{0} \delta\bigl(x - (x_j(t) + \lambda \, (x_i(t) - x_j(t))) \bigr) \mathrm{d} \lambda.
\end{align*}
To get averaged quantities that can be passed to the macroscopic model, we have to average the distributions $\widehat{\rho}$, $\widehat{v}$ and $\widehat{p}$ over a sampling domain. Consequently we obtain the spatially averaged, microscopic quantities $\rho$, $v$, and $p$.
In the following we will only consider these averaged quantities.

For a homogeneous particle chain with constant particle masses $m = m_i$, the averaged microscopic pressure is given by $p(\tau) = - \phi'(\tau)$, as a function of the specific volume $\tau = m/\rho$, if the local microscopic temperature is zero, which is the case in our setting, as the particles are initialized without any random fluctuations.
Using this relation, the macroscopic pressure function can be determined directly from the microscale model. 
This means that for the consistency of both models we have to set $\phi(\tau) = \psi(\tau)$, where $\psi$ denotes the specific Helmholtz free energy of the macroscopic system, satisfying $p(\tau) = - \psi'(\tau)$.
In the following, we consider the potential 
\begin{align}  \label{eq:vdw-potential}
  \phi(r)  &= -\frac{a}{r} - R\theta \, \operatorname{ln}(b-r), 
  &
  \phi'(r) &= \frac{a}{r^2} + \frac{R\theta}{b-r}, 
\end{align}
which is consistent with the van der Waals pressure \eqref{eq:van-der-waals-pressure}. 
However, we stress that the choice of the potential is arbitrary and implies the macroscopic pressure, not the other way round. 
Here, the explicit choice of $\phi$ is done to compare the multiscale scheme with already existing solvers for van der Waals fluids.

\subsection{The Microscopic Riemann Problem}
\label{sec:microscopic-riemann-problem}

\begin{figure}[t]
\sidecaption
\centering
\mbox{}\hfill
\includegraphics[width=0.5\columnwidth]{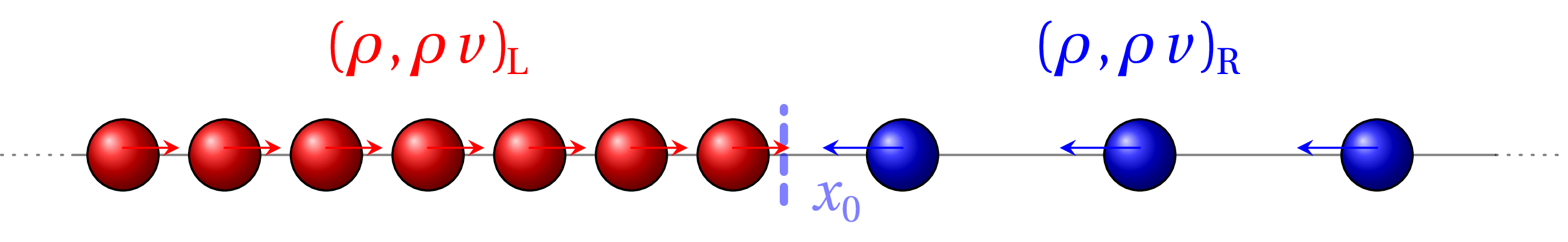}
\hfill\mbox{}
\caption{Schematic representation of Riemann data at the microscopic scale.}
\label{fig:particlechain_riemann}
\end{figure}

Our main goal is to describe the dynamics of the fluid at the liquid-vapor interface, which we interpret as a Riemann problem. 
To incorporate microscopic properties, we define a Riemann problem on the microscopic scale and solve it in order to extract the wave pattern, which will be used to compute the fluxes at the interface on the macroscopic scale.
\newline
Therefore we have to convert the macroscopic quantities to the microscale quantities and vice versa using the Irving--Kirkwood formulas \eqref{eq:irving-kirkwood:1}.
To be more precise, for macroscopic Riemann problem data $u_{\mathrm{L}} = (\rho, \rho v)_{\mathrm{L}}$ and $u_{\mathrm{R}} = (\rho, \rho v)_{\mathrm{R}}$ we set the initial particle configuration uniformly, such that for both $\alpha = \mathrm{L}$ and $\alpha = \mathrm{R}$
\begin{align*}
 x^0_{i} - x^0_{i-1} &= \frac{m_i}{\rho_{\alpha}},  
 &
 v^0_i &= v_{\alpha}, &
 &\text{for all } i \in I_{\alpha},
\end{align*}
holds, where $I_{\mathrm{L}} = \{i ~\vert~ i = 1,\ldots,N ~ \text{ with } x_i \leq 0\}$, $I_{\mathrm{R}} = \{i ~\vert~ i = 1,\ldots,N ~ \text{ with } x_i > 0\}$ are the index sets for the left/right hand particles.
A schematic depiction of such a configuration can be seen in Figure \ref{fig:particlechain_riemann}.
This gives us the microscopic Riemann problem for
\begin{align}
 ({\rho}, {\rho}{v})(x, \, t=0) 
 = 
 \begin{cases}
  ({\rho}, {\rho}{v})_\mathrm{L}
  &: ~ x \leq 0, \\
  ({\rho}, {\rho}{v})_\mathrm{R}
  &: ~ x > 0, 
 \end{cases}
\end{align}
with a left state $({\rho}, {\rho}{v})_\mathrm{L}$ and a right state $({\rho}, {\rho}{v})_\mathrm{R}$, defined by local averages of \eqref{eq:irving-kirkwood:1}, with the jump at zero. 
\newline
After running the microscale simulation, the evolving wave pattern has to be transferred to the macroscopic model.
For that we perform some local averaging over the particles states using the Irving--Kirkwood formulas \eqref{eq:irving-kirkwood:1}.
The interface speed is obtained by tracking the interface position on the microscopic scale.

\subsubsection{Extracting Key Quantities}
\label{sec:extracting_key_quantities}

For a given solution to the microscopic Riemann problem, it is still the question how we extract the key quantities from the microscopic solution. 

\begin{figure}
\sidecaption
\centering
\mbox{}\hfill
\includegraphics[width=0.6\columnwidth]{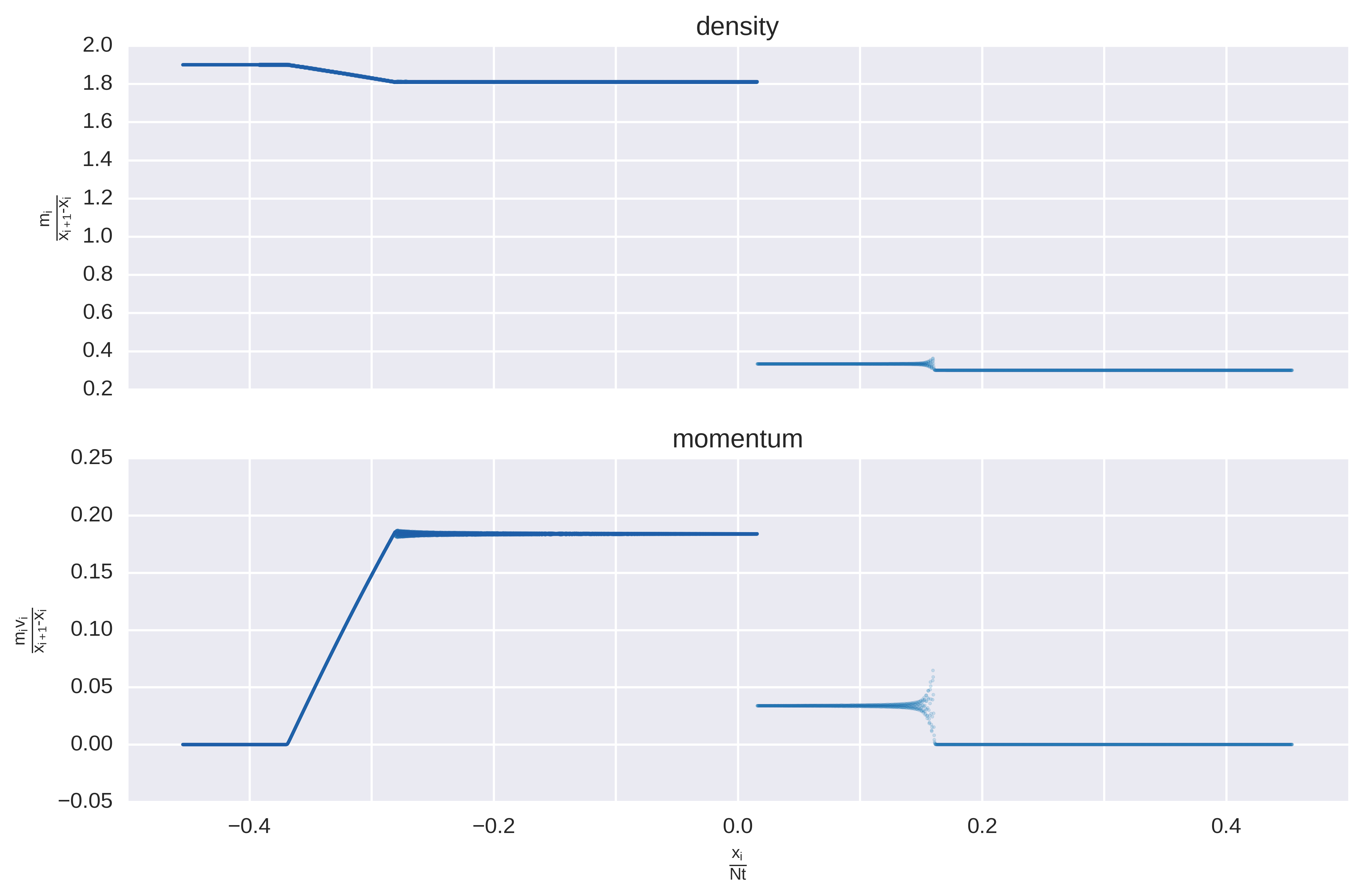}
\hfill\mbox{}
\caption{Example of a solution of the microscopic particle model with van der Waals potential \eqref{eq:vdw-potential} for the initial values $({\rho}, {\rho}{v})_\mathrm{L} = (1.9, 0)$  and $({\rho}, {\rho}{v})_\mathrm{R} = (0.3, 0)$ for $16000$ particles at $t = 2500$. 
The phase boundary is located at the density jump near the origin.}
\label{fig:particlesolution}
\end{figure}

In Figure \ref{fig:particlesolution} an example of a solution of the particle model is depicted.
It can be seen that, similar to wave patterns in the continuum case, a 3-wave pattern evolves -- see Figure \ref{fig:wavepattern}. 
We apply this analogy to construct a numerical flux for the interface dynamics. 
To this end, similar to the numerical flux in \cite{chalons.rohde.ea:finite:2016}, we need to extract the states adjacent to the interface from the wave pattern, and also the interface propagation speed.  
To obtain these values, the interface is tracked by considering the biggest local change in density and then the neighboring states can be computed easily by local averaging left and right of the interface.

\subsection{Discretization of the Particle System}
For the time-discretization of the particle system we apply the velocity Verlet algorithm \cite{verlet:computer:1967}. 
It is an explicit scheme with microscale time step $\Delta t > 0$ of the following form:
\begin{align}
\label{eq:vel-verlet}
\begin{aligned}
 x(t+\Delta t) &= x(t) + \Delta t \, v(t)  + \tfrac{1}{2} \Delta t^2 \, a(t), \\
 v(t+\Delta t) &= v(t) + \tfrac{1}{2} \Delta t \, \bigl( a(t) + a(t + \Delta t) \bigr),
\end{aligned}
\end{align}
where $v = \tfrac{\mathrm{d} x}{\mathrm{d} t}$ is the particle velocity, and $a = \tfrac{\mathrm{d} v}{\mathrm{d} t}$ the particle acceleration, computed from the forces between the particles at each time step.
It is of second order and has the advantage that no intermediate values of $x$, $v$, or $a$ have to be stored.
Furthermore, we see that all steps can be run in parallel. 
This enables us to run the particle simulations on a graphics processing unit (GPU) which gives a major speedup, as opposed to conventional hardware.

\section{The Multiscale Model}
To design the multiscale model, we consider the continuum model \eqref{eq:euler-2phase} with the interface conditions \eqref{eq:interface_conditions} as our macroscopic model.  
The bulk phases of the continuum model are solved by a standard finite volume scheme, and we focus on the description of the interface dynamics. 
We refrain from formulating a kinetic relation, and instead include data from the microscopic Riemann solutions of the particle model presented in Section \ref{sec:microscopic-riemann-problem}.  
Hereby, the communication between the macroscale continuum model and microscale particle model is solely data driven.
Only the macroscopic constraints $({\rho}, {\rho}{v})_\mathrm{L}$ and  $({\rho}, {\rho}{v})_\mathrm{R}$ are needed for setting up the microscale Riemann problem, and in return, for the computation of the macroscale interface flux just the response values $(s, u^{*}_{\mathrm{L}}, u^{*}_{\mathrm{R}})$ from the wave pattern are needed, see Section \ref{sec:discretization-continuum}. 
Consequently, for the continuum model only the input-output relation $(u_{\mathrm{L}}, u_{\mathrm{R}}) \mapsto (s, u^{*}_{\mathrm{L}}, u^{*}_{\mathrm{R}})$ from the microscopic Riemann problem is important.

\subsection{Model Reduction Algorithm}
\label{sec:model_reduction}
 
The evaluation of the microscale model is computationally relatively expensive, and if it is evaluated at each interface edge and time step of the continuum model, the coupled micro- macroscale model becomes computationally unfeasible - see Section \ref{sec:numerical_simulations} for a more details. 
To counter this problem we exploit the fact that the coupling is solely data-driven, and apply a reduced, kernel-based surrogate model for the particle model input-response relation 
$f_{\mathrm{micro}} \colon
 (u_{\mathrm{L}}, u_{\mathrm{R}}) \mapsto (s, u^{*}_{\mathrm{L}}, u^{*}_{\mathrm{R}})$,
where $u := (\rho, \rho v)$.
More abstractly, we apply the microscale model as a black box and put the reduced model into the framework of machine learning. 
For that $x \in \mathbb{R}^{d_1}$ denotes the $d_1$-dimensional input data, which is in our case $x = (u_{\mathrm{L}}, u_{\mathrm{R}})$, and $y \in \mathbb{R}^{d_2}$ is the $d_2$-dimensional response of our model, in our case the measured data $(s, u^{*}_{\mathrm{L}}, u^{*}_{\mathrm{R}})$. 
The aim is now, to train a regression function from samples a set $D_n = \{(x_i,y_i) \, : \, i = 1, \ldots, n\}$, obtained from observations
$y_i = f_{\mathrm{micro}}(x_i) + \varepsilon_{\mathrm{s}}$ 
that describes $f_{\mathrm{micro}}$ in an optimal sense.
Here $\varepsilon_{\mathrm{s}}$ accounts for possible normal distributed measurement noise. 
To get the regression function from the sample set $D_n$ we apply a support vector regression scheme, see e.g. \cite{steinwart.christmann:support:2008}.
Therefore we have to train the reduced model function 
\begin{align*}
 f(x) = \sum_{i=1}^{n} \alpha_i \, k_\gamma(x_i, x), 
\end{align*}
on the trainings data set $D_n$, where $k_\gamma$ is the radial basis kernel function $k_\gamma(x_i,x) = \exp(-\gamma \|x-x_i \|^2)$
In this context, that means that we have to determine the coefficients $\alpha_i \in \mathbb{R}$ such that $f$ describes $f_{\mathrm{micro}}$ optimally under the observations in $D_n$. 
Consequently, an optimization problem has to be solved each time the reduced model is trained. 
\newline
More details on kernel-based surrogate modeling can be found in e.g. \cite{kissling.rohde:computation:2015,wirtz.karajan.ea:surrogate:2015}.

\subsubsection{Dynamic Sampling Scheme}

\begin{figure}[h]
\centering
\includegraphics[width=0.75\columnwidth]{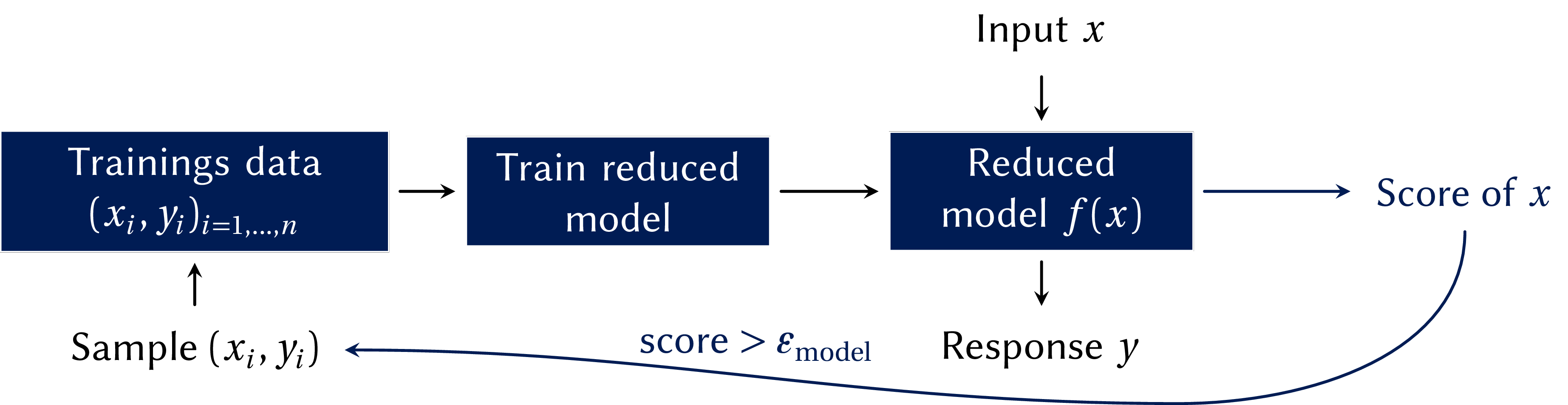}
\caption{Sketch of the model reduction scheme with dynamic sampling.}
\label{fig:model_reduction_scheme}
\end{figure}

In our case, the input values that are needed cannot be prescribed a priori. 
Therefore we apply a dynamic sampling strategy, which is described in this section.
\newline
The sampling set $D_n$ is updated dynamically at each time step of the continuum model. 
To this end, we assign each input value $x \in \mathbb{R}^{d_1}$ a score $\gamma(x; D_n)$ that describes the quality of the surrogate model at the point $x$. 
This score is computed at each evaluation of the surrogate model.
If the score is below a certain threshold $\varepsilon_{\mathrm{model}} > 0$, we simply evaluate the point $x$ by the surrogate model. On the other hand, if it is above the threshold, we draw a new sample by evaluating the microscale model and add it to the training set $D_{n+1} = D_n \cup \{(x_{n+1},y_{n+1})\}$. 
A sketch of the complete model reduction scheme is shown in Figure \ref{fig:model_reduction_scheme}.
\newline
In the following we use the distance from an input value $x$ to the nearest point of the sample data set $D_n$, i.e.
$ \gamma(x\,;\, D_n) = \min_{i \leq n} \|x - x_i\|$. 
One drawback of this simple choice is, that we only consider the input values and ignore the output values, which could give an indication whether the (local) variance of the underlying model equation is higher or lower in certain areas of the input space.   

\subsection{Numerical Discretization of the Multiscale Model}
\label{sec:discretization-continuum}

To discretize the macroscale model, we apply the time-explicit front tracking finite volume scheme for systems from \cite{chalons.rohde.ea:finite:2016}. 
It has the advantage that the sharp interface is resolved within the mesh, i.e. the discretized phase boundary always coincides with a (moving) mesh edge.
At the interface we have to solve a special Riemann problem including the phase dynamics. 
From its solution we have to extract the interface propagation speed $s$ and the adjacent fluid states $u^{*}_{\mathrm{R}}$ and $u^{*}_{\mathrm{L}}$, see Figure \ref{fig:wavepattern}.
However, instead of solving the microscale Riemann problem each time, we insert the model reduction scheme from Section \ref{sec:model_reduction}. 
The wave pattern values are inserted in the numerical flux at the interface $g(u_\mathrm{L},u_\mathrm{R}) = \tfrac{1}{2} \left( f(u^{*}_{\mathrm{L}}) + f(u^{*}_{\mathrm{R}}) - s(u^{*}_{\mathrm{L}} + u^{*}_{\mathrm{R}}) \right)$. 
In the bulk phases we apply a standard Lax--Friedrichs flux scheme.

\section{Numerical Simulations}
\label{sec:numerical_simulations}

In this section we present some numerical simulation results to show that the multiscale scheme is viable and applicable to (two-dimensional) droplet dynamics.

\begin{figure}
\sidecaption
\centering
\mbox{}\hfill
\includegraphics[width=0.6\columnwidth]{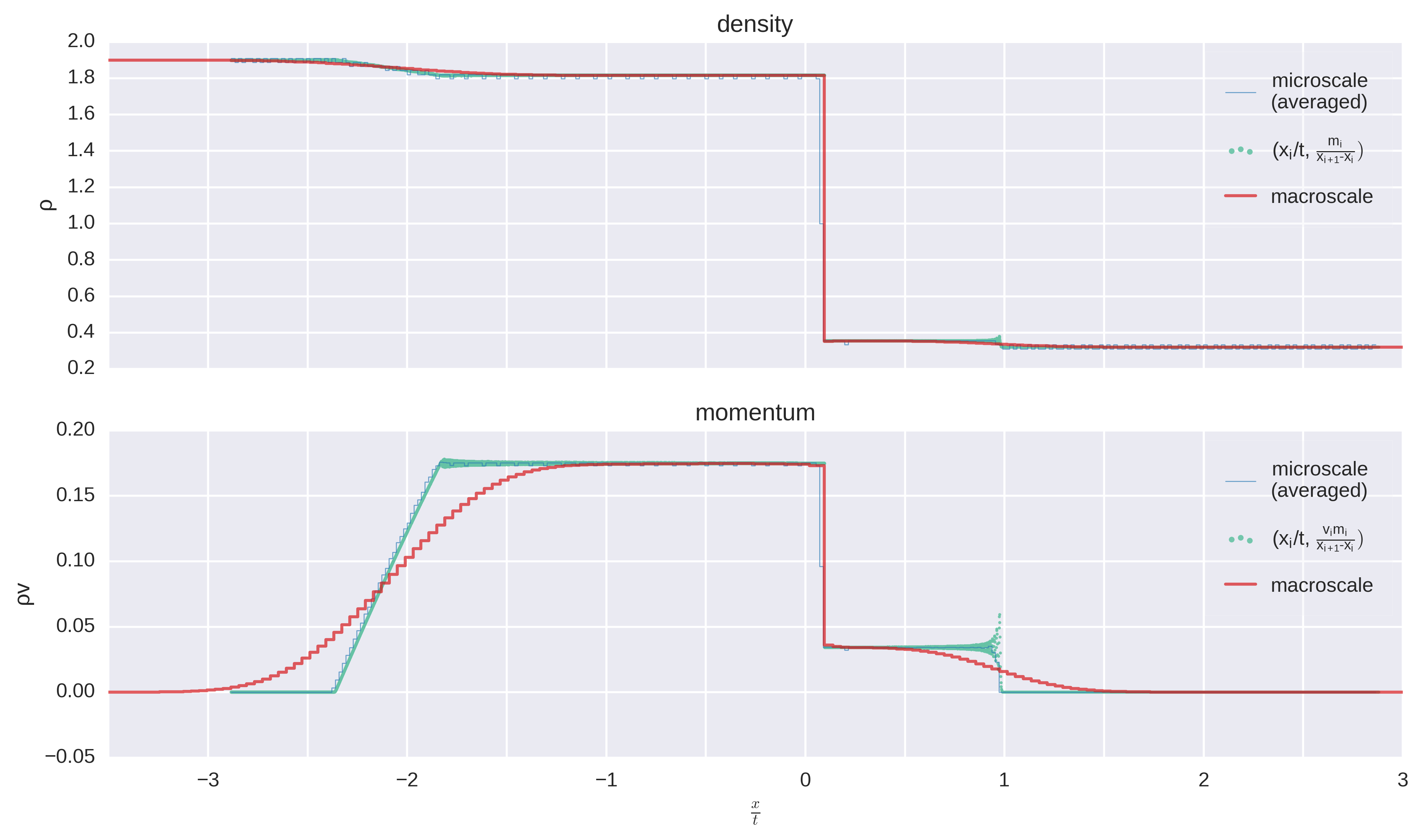}
\hfill\mbox{}
\caption{One-dimensional solution of the multiscale model and the particle model for the Riemann problem, where the phase boundary is located at the density jump near $x/t = 0.2$.}
\label{fig:one_dim_vdw_euler}
\end{figure}

\paragraph{A Multiscale Simulation of the Riemann Problem:} 
The first simulation results show the consistency between the particle model and the multiscale model in one spatial dimension.
Therefore, we run both, the particle model and the multiscale model for the same set of Riemann data and compare the averaged particle solution with the multiscale solution. 
For the initial conditions we have $\rho_{\mathrm{L}} = 2.0$, $v_{\mathrm{L}} = 0$ for $x < 0$ in the liquid phase, and on the right side the vapor-phase Maxwell equilibrium state $\rho_{\mathrm{R}} \approx 0.317$, $v_{\mathrm{R}} = 0$.
In Figure \ref{fig:one_dim_vdw_euler} 
both solutions are superimposed and we can see that they fit well, and in particular the wave speeds of the phase boundary coincide.

\paragraph{Multiscale Simulations of a Droplet in 2D:}
Next, we solve the multiscale model on the continuum scale in two spatial dimensions. 

\begin{figure}[h]
\centering
\hfill%
\includegraphics[width=0.32\columnwidth]{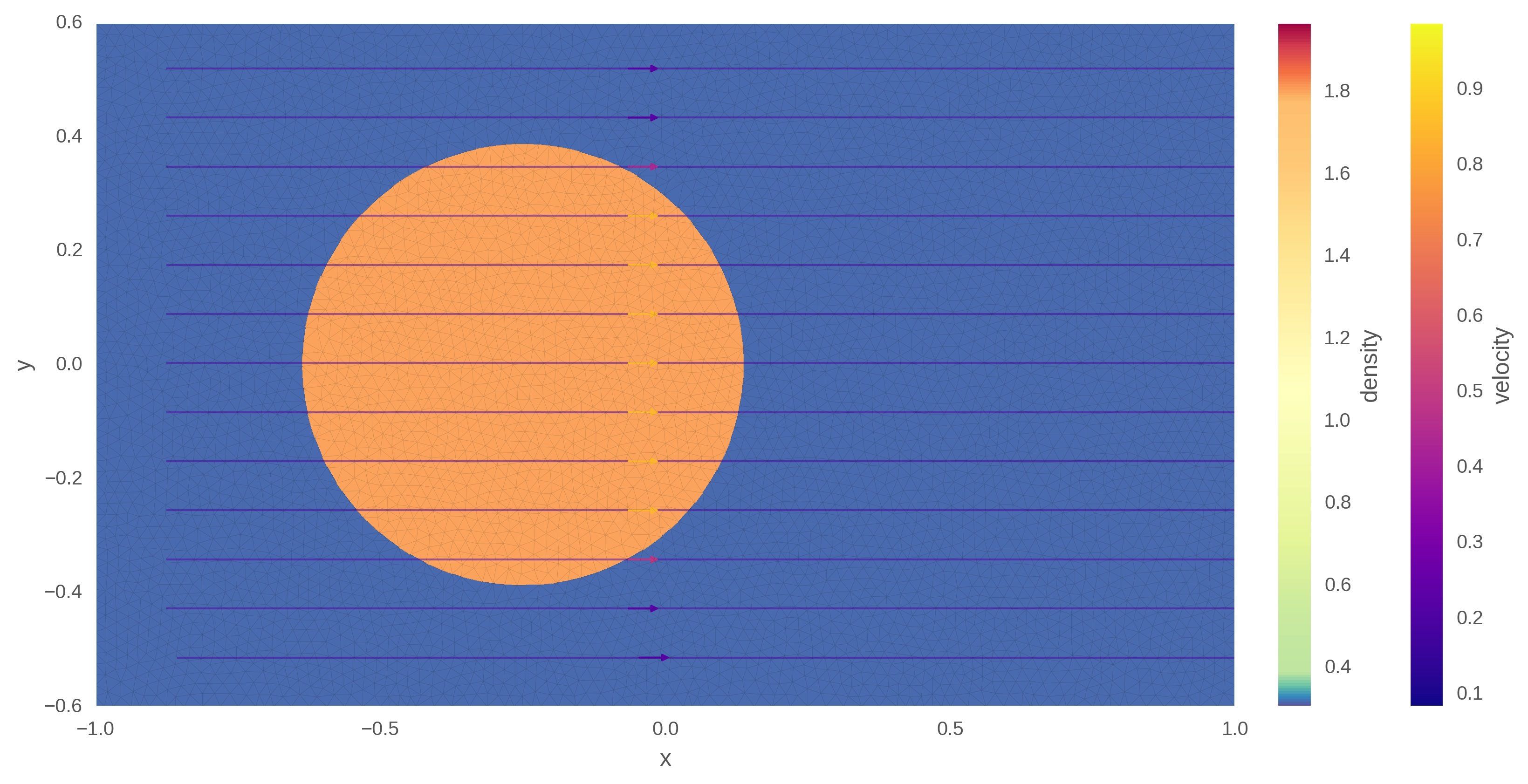}%
\hfill%
\includegraphics[width=0.32\columnwidth]{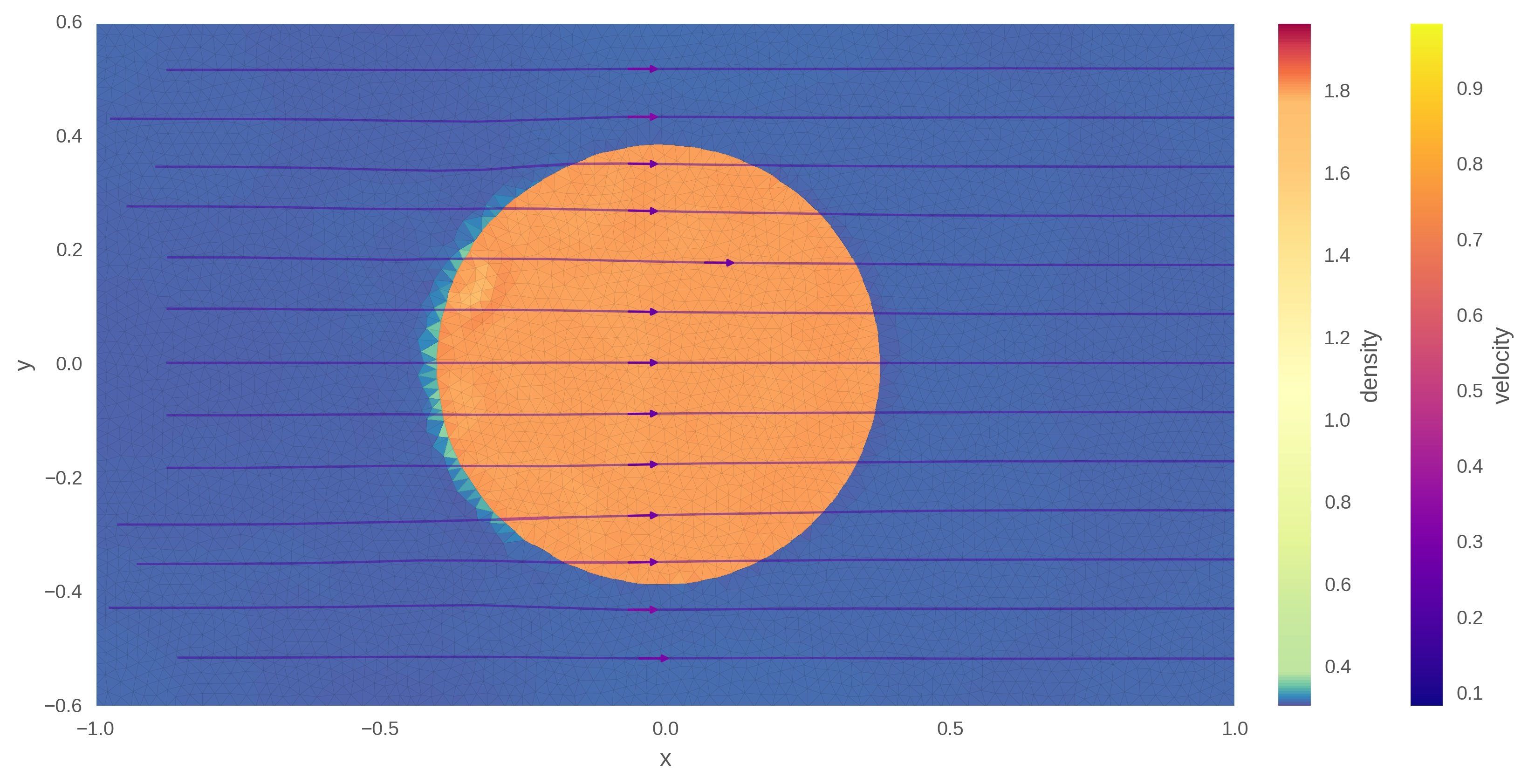}%
\hfill%
\includegraphics[width=0.32\columnwidth]{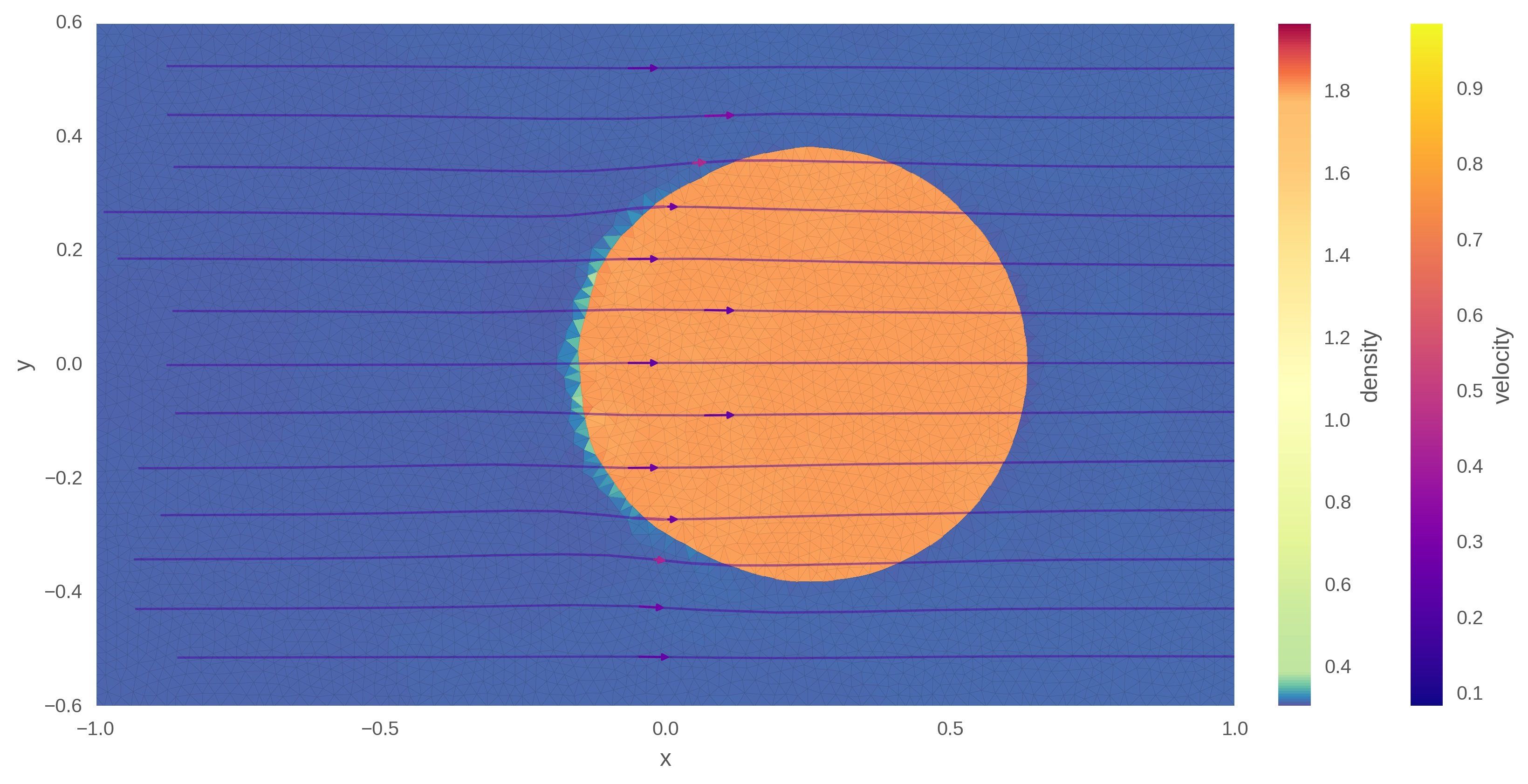}%
\hfill%
\caption{Multiscale simulation of a moving droplet at $t = 0$, $t = 1.25$ and $t = 2.5$ (from left to right).}
\label{fig:moving_droplet}
\end{figure}

\textit{Droplet transport:}
In the first simulation we present the performance of the front tracking scheme in two spatial dimensions. 
The initial conditions for the density are the Maxwell equilibrium states, which are $\rho_{\mathrm{liq}} \approx 1.804$ for the liquid phase and $\rho_{\mathrm{vap}} \approx 0.317$ for the vapor phase. 
The initial velocity in the domain and on the boundary is set to $v = (0.2, 0)^\top$. 
In Figure \ref{fig:moving_droplet} we see that the droplet is transported through the domain and mostly keeps its shape.
Furthermore it remains in equilibrium, the increased density at th interface in the vapor on the left side and the small oscillations are due to the local averaging if the triangulation is restructured. 

\begin{figure}
\centering
%\resizebox{0.96\columnwidth}{!}{
\hfill%
\includegraphics[height=2.8cm]{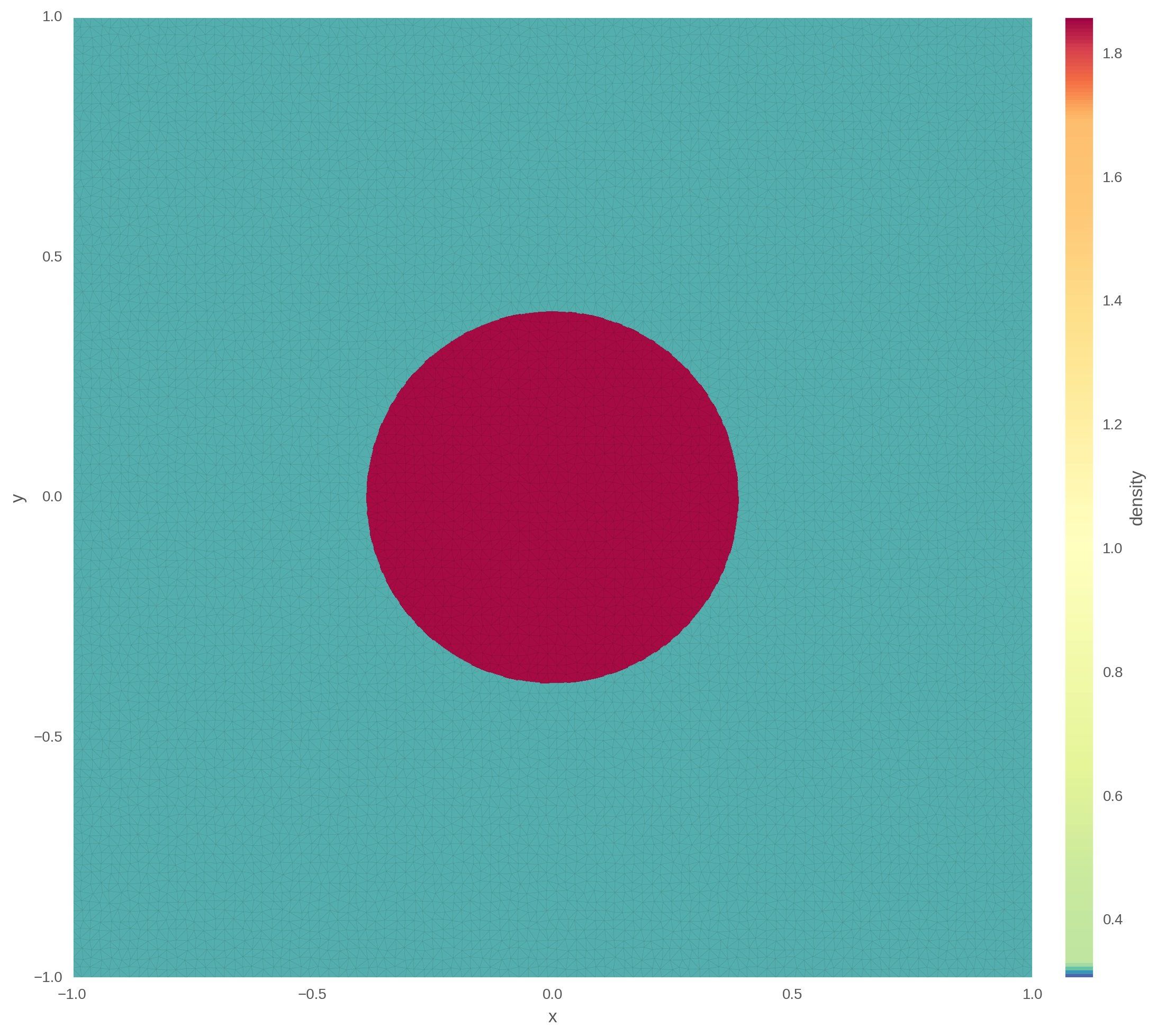}%
\hfill%
\includegraphics[height=2.8cm]{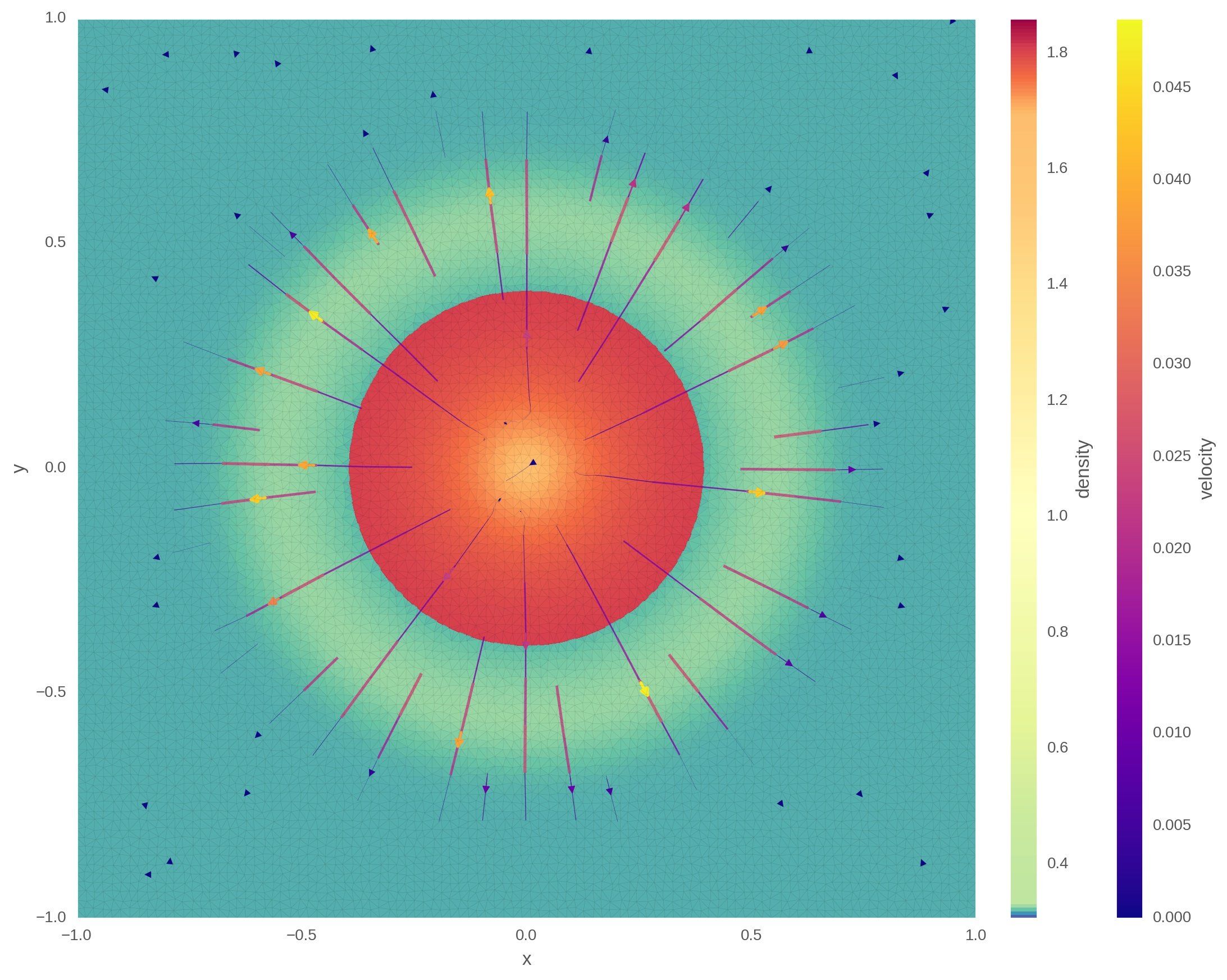}%
\hfill%
\includegraphics[height=2.8cm]{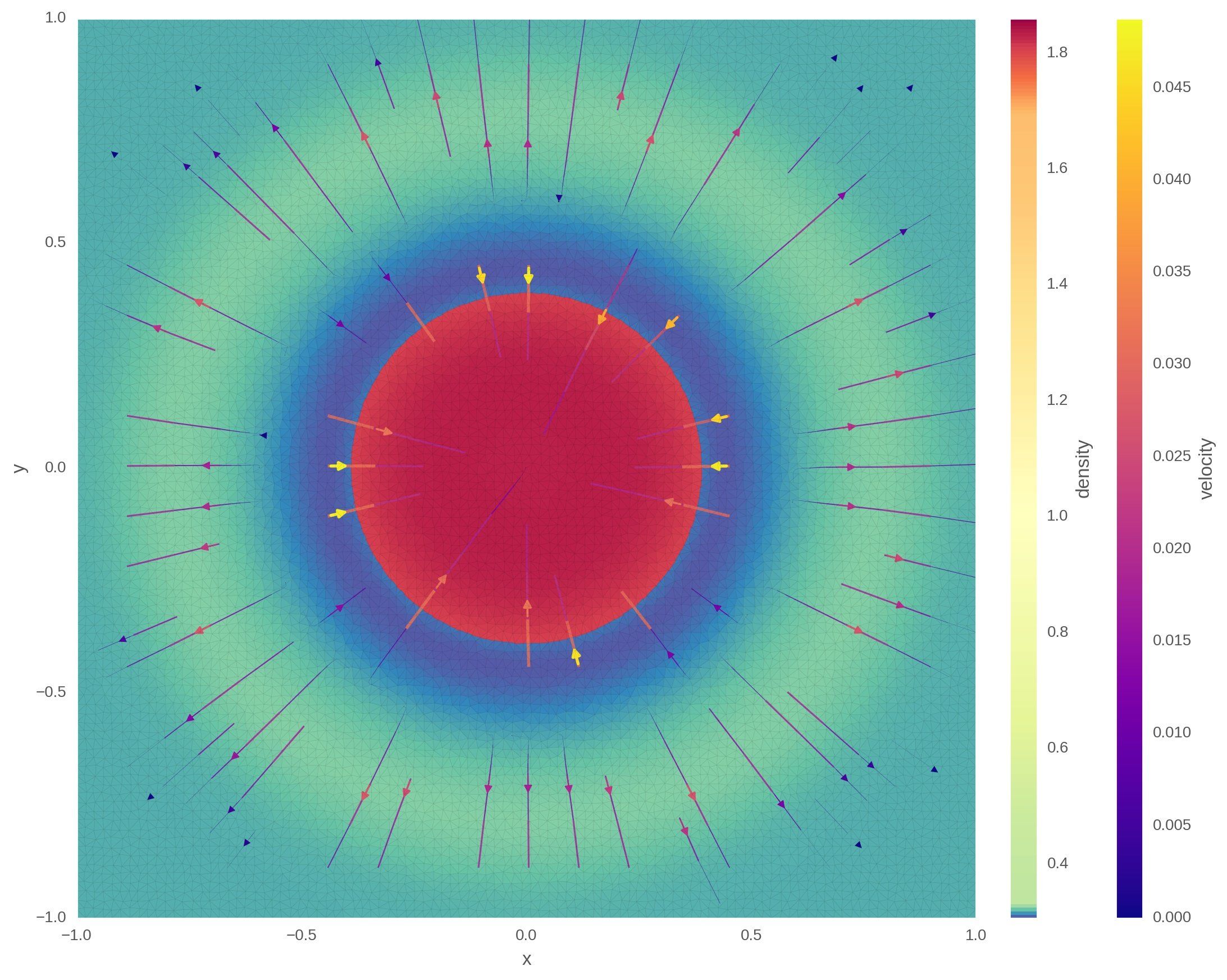}%
\hfill%
\includegraphics[height=2.8cm]{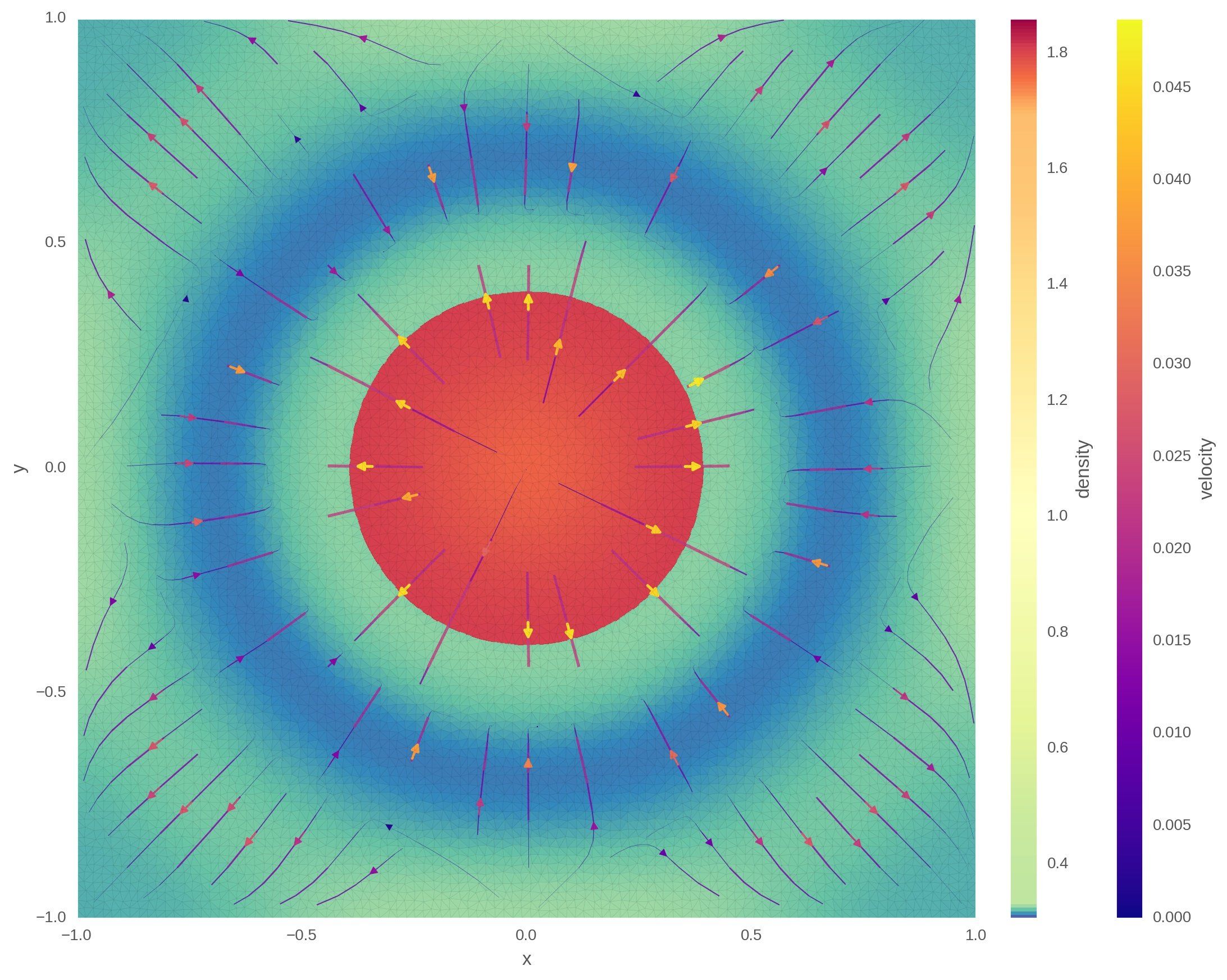}%
\hfill%
%}
\caption{Multiscale simulation of an oscillating droplet at $t = 0, 0.25, 0.5, 0.75$ (from left to right).}
\label{fig:oscillating_droplet}
\end{figure}

\begin{figure}
% \sidecaption
\centering
\includegraphics[width=0.45\columnwidth]{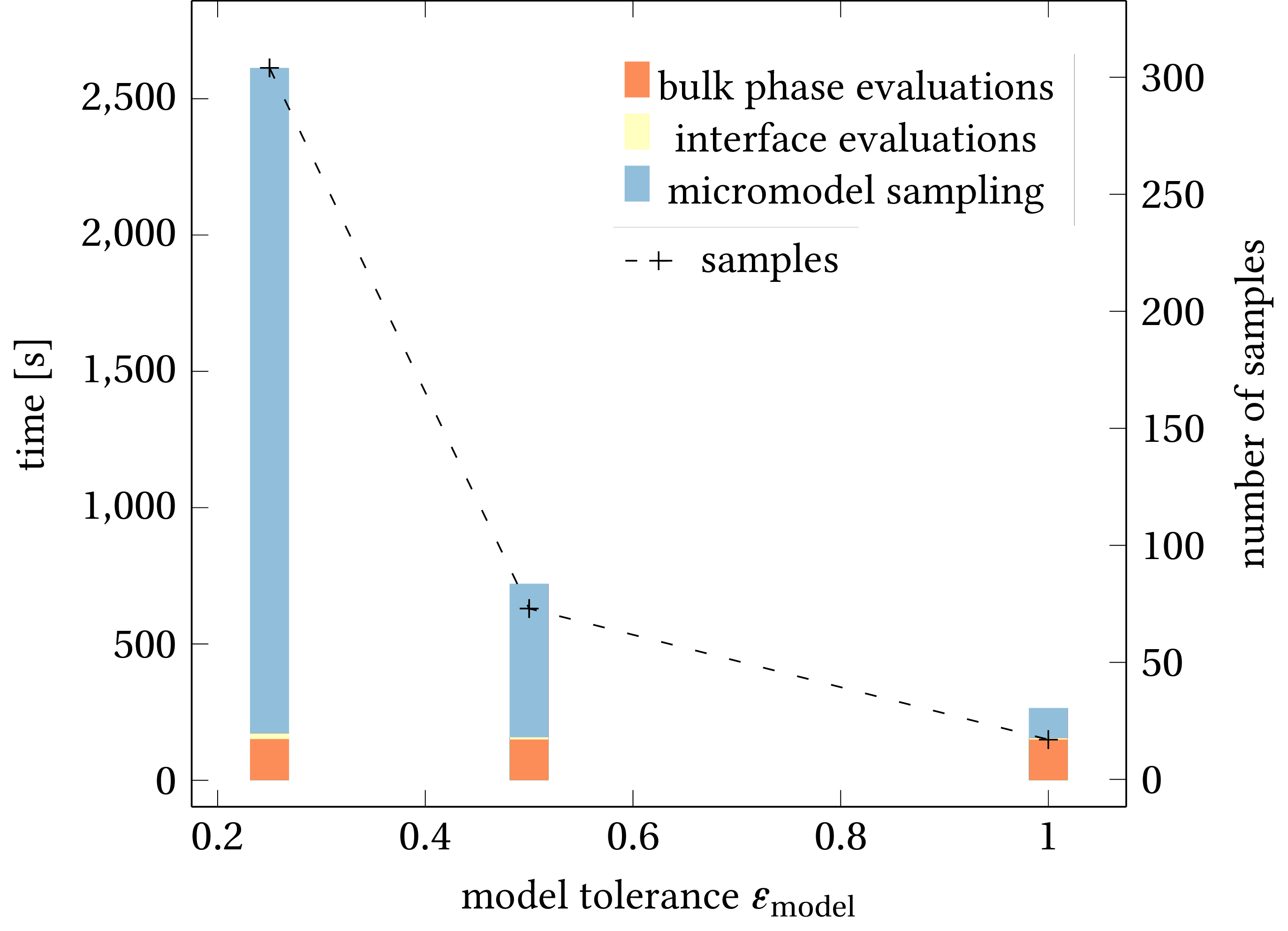}%
\caption{Computational time in seconds with respect to the model tolerance $\varepsilon_{\mathrm{model}}$. The dashed line indicates the number of samples that are drawn from the microscale model.}
\label{fig:computational_time} 
\end{figure}

\textit{Oscillating droplet:}
In the next simulation we consider a droplet that is perturbed from the liquid phase equilibrium, i.e. $\rho_{\mathrm{liq}} = 1.85$, and measure the effect of the model tolerance $\varepsilon_{\mathrm{model}}$ on the computational time. 
The simulation results for $\varepsilon_{\mathrm{model}} = 0.5$ are presented in Figure \ref{fig:oscillating_droplet}.
We consider reflecting boundary conditions and thus, the droplet oscillates slightly. 
The computational time\footnote{All simulations were performed on a single workstation equipped with an Intel\textsuperscript{\textregistered} i7-6700 CPU at 3.4 GHz, 16GB RAM, and a Nvidia\textsuperscript{\textregistered} GTX980 Ti GPU.} for this simulation is depicted in Figure  \ref{fig:computational_time}.
It can be seen that the time for computing new samples is of the same order of the finite volume computations, which underlines the performance of the model reduction scheme. 
If we would not apply model reduction, we would have to run microscale simulation (around 20 seconds per sample) for all 8000 time steps at each of the $\sim$160 interface edges. 
This would lead to a computational time that amounts to roughly one year. 
Compared to that the runtime with the model the model reduction scheme takes only several minutes. 
This gives us huge speedups (with/without model reduction) as shown in Table \ref{tab:speedup}.

\begin{table}[h]
\centering
\begin{tabular}{lr}
\toprule
 $\varepsilon_{\mathrm{model}}$ & speedup \\
 \midrule
1.0 & 222112 \\
0.5 &  47000 \\
0.25 & 10836 \\ 
\bottomrule
\end{tabular}
\caption{Speedup, with/without model reduction.} 
\label{tab:speedup}
\end{table}

\section{Conclusions}
In this work, we have presented a multiscale model for the description of two-phase flows with a sharp interface, that incorporates microscale features originating from an atomistic particle model. 
We have exploited the fact that the coupling of the micro- and macroscale model is solely data-based and developed a model reduction scheme that dynamically draws new data points from the microscale model and makes the whole multiscale scheme computationally feasible. 
Numerical simulation results are presented, that not only showing the consistency of the multiscale scheme, but also that the applicability in more complex situations without prescribing some (ad-hoc) kinetic relations. 

\begin{acknowledgement}
The work was supported by the German Research Foundation (DFG) through SFB TRR 75 ``Droplet dynamics under extreme ambient conditions''. 
\end{acknowledgement}

\bibliographystyle{plain}
\bibliography{bibliography}

\end{document}